\documentclass[12pt]{article}

\usepackage{amsmath}
\usepackage{amsfonts}
\usepackage{amssymb}
\usepackage{amsthm}
\usepackage{authblk}
\usepackage{bm}
\usepackage{booktabs} 
\usepackage[font=scriptsize]{caption}
\usepackage{color}
\usepackage{esvect}
\usepackage{enumitem}
\usepackage{extarrows}
\usepackage{float}
\usepackage{footmisc} 
\usepackage{graphicx}
\usepackage{hhline}
\usepackage[utf8]{inputenc}
\usepackage{hyperref} 
\usepackage{indentfirst}
\usepackage{mathrsfs}
\usepackage{multirow}
\usepackage{modroman}
\usepackage{natbib}
\usepackage{geometry}
\usepackage{pdfpages}
\usepackage{setspace}
\usepackage{subcaption}
\usepackage{tabularx}
\usepackage{url}
\usepackage{xcolor}
\usepackage{xspace}

\newtheorem{theorem}{Theorem}

\newtheorem{condition}{Condition}
\newtheorem{corollary}{Corollary}
\newtheorem{remark}{Remark}

\makeatletter
\newcommand{\distas}[1]{\mathbin{\overset{#1}{\kern\z@\sim}}}%
\newsavebox{\mybox}\newsavebox{\mysim}
\newcommand{\distras}[1]{%
  \savebox{\mybox}{\hbox{\kern3pt$\scriptstyle#1$\kern3pt}}%
  \savebox{\mysim}{\hbox{$\sim$}}%
  \mathbin{\overset{#1}{\kern\z@\resizebox{\wd\mybox}{\ht\mysim}{$\sim$}}}%
}
\makeatother
\newcommand\numberthis{\addtocounter{equation}{1}\tag{\theequation}}

\providecommand{\keywords}[1]
{
  \small	
  \textbf{\textit{Keywords:}} #1
}

\numberwithin{equation}{section}
\numberwithin{theorem}{section}
\numberwithin{lemma}{section}
\numberwithin{prop}{section}
\numberwithin{condition}{section}
\numberwithin{corollary}{section}
\numberwithin{remark}{section}
\numberwithin{example}{section}
\numberwithin{definition}{section}

\title{Asymptotic Normality of the Conditional Value-at-Risk based Pickands Estimator}
\author[a]{\normalsize Yizhou Li\thanks{Corresponding author at Department of Applied Mathematics and Statistics at Stony Brook University,
United States. E-mail address: \href{mailto:yizhou.li@stonybrook.edu}{yizhou.li@stonybrook.edu}}}
\author[a,b]{\normalsize Pawe\l\ Polak}

\affil[a]{\small \textit{Department of Applied Mathematics and Statistics, Stony Brook University, United States}}
\affil[b]{\small \textit{Institute for Advanced Computational Science, Stony Brook University, United States}}

\setcitestyle{authoryear,open={(},close={)}} 
\begin{document}
\maketitle

\begin{abstract}
The Pickands estimator for the extreme value index is beneficial due to its universal consistency, location, and scale invariance, which sets it apart from other types of estimators. However, similar to many extreme value index estimators, it is marked by poor asymptotic efficiency. \citet{teng2021thesis} introduces a Conditional Value-at-Risk (CVaR)-based Pickands estimator, establishes its consistency, and demonstrates through simulations that this estimator significantly reduces mean squared error while preserving its location and scale invariance. The initial focus of this paper is on demonstrating the weak convergence of the empirical CVaR in functional space. Subsequently, based on the established weak convergence, the paper presents the asymptotic normality of the CVaR-based Pickands estimator. It further supports these theoretical findings with empirical evidence obtained through simulation studies.
\end{abstract}

\section{Introduction} \label{section:intro}
Suppose $X$ is a random variable (r.v.) with a distribution function (d.f.) $F$, and $X_{1}, X_{2}, \ldots, X_{n}$ is a sequence of independent samples of $X$. The extreme value d.f. with shape parameter $\gamma\in\mathbb{R}$ (shift parameter $0$ and scale parameter $1$) is defined by
\begin{align*}
    G_{\gamma}(x) = \begin{cases}
        \exp \left\{-(1+\gamma x)^{-1 / \gamma}\right\}, & \quad \text{for } \gamma\neq0,\ 1+\gamma x>0, \\
        \exp \left\{-e^{-x}\right\}, & \quad \text{for } \gamma=0,\ x\in\mathbb{R}.
    \end{cases}
\end{align*}
The d.f. $F$ belongs to the \textit{domain of attraction} of $G_\gamma$ for some $\gamma\in \mathbb{R}$, noted as $F\in\mathscr{D}(G_\gamma)$, if there exists a sequence of constants $a_n>0$ and $b_n \in \mathbb{R}$ such that
\begin{equation} \label{eqn: DoA}
    \lim _{n \rightarrow \infty} P\left\{\max \left(X_{1}, \ldots, X_{n}\right) \leqslant a_{n}x + b_{n}\right\}=\lim _{n \rightarrow \infty} F^{n}\left(a_{n} x+b_{n}\right)=G_{\gamma}(x).
\end{equation}
This $\gamma$ is called the \textit{extreme value index} of $F$.

An alternative characterization of $\gamma$ in extreme value theory can be provided using the generalized Pareto distribution (GPD), which has d.f.
\begin{align*}
    H_\gamma(x)=\begin{cases}
        1-(1+\gamma x)^{-1 / \gamma}, & \quad \text{for }\gamma\neq0, \ x\geq0, \ 1+\gamma x>0, \\
        1 - e^{-x}, & \quad \text{for } \gamma=0, \ x\geq0.
    \end{cases}
\end{align*}
Note that $H_\gamma$ with $\gamma>0$ yields a Pareto distribution; $H_0$ coincides with the exponential distribution. Additionally, $H_\gamma$ with $\gamma<0$ exhibits a finite right endpoint, particularly with $\gamma=-1$ resulting in a uniform distribution on the interval $(0,1)$.

Understanding the extreme value index $\gamma$ is crucial for modeling maxima and estimating extreme quantiles. The estimation of $\gamma$, besides high quantile estimation, is one of the most crucial problems in univariate extreme value theory.

Let $X_1^{(n)},X_2^{(n)},\cdots,X_n^{(n)}$ denote the descending order statistics of $X_1,\ldots,X_n$. \cite{pickands1975statistical} proposed a simple estimator for $\gamma$, constructed in terms of log-spacings of order statistics:
\[
    \widehat{\gamma}_{n}^{P}(m):=\frac{1}{\log 2} \log \left(\frac{X_{m}^{(n)}-X_{2 m}^{(n)}}{X_{2 m}^{(n)}-X_{4 m}^{(n)}}\right), \quad \text{for } m=1,\ldots,\lfloor n/4 \rfloor,
\]
where $\lfloor x \rfloor$ is the floor of a real number $x$. Compared with the \cite{hill1975simple} estimator and the probability weighted moment estimator proposed by \cite{hosking1985estimation}, the primary benefits of the Pickands estimator include the invariant property under location and scale shift, and the consistency for any $\gamma \in \mathbb{R}$ under intermediate sequence such that $m_{n} \to \infty$, $m_{n} / n \to 0$ as $n\to\infty$.

However, the Pickands estimator exhibits relatively poor asymptotic variance. To address the drawback, 
\cite{yun2002generalized} first generalized the Pickands estimator by
\begin{equation} \label{eqn: yun}
    \widehat{\gamma}_{n, m}^{Y}(u, v):=\frac{1}{\log v} \log \frac{X_{m}^{(n)}-X_{[u m]}^{(n)}}{X_{[v m]}^{(n)}-X_{[u v m]}^{(n)}}, \quad u, v>0, u, v \neq 1,
\end{equation}
where $1 \leqslant m,[u m],[v m],[u v m] \leqslant n$ and $[x]$ denotes the integer part of $x \in \mathbb{R}$. The estimator is constructed by the linear combinations of the logarithm of the spacing between intermediate order statistics. The spacings are generalized and determined by $u$ and $v$. The optimal $(u^\ast, v^\ast)$ with $0<v\leq u<1$ along $\gamma$ which minimizes the asymptotic variance is approximated numerically. Then, the estimation is conducted in an adaptive manner. 

All the previously mentioned estimators are constructed using order statistics as empirical quantiles of the distribution function (d.f.). Alternatively, a different quantity called Conditional Value at Risk (CVaR) order statistic, also known as the empirical super-quantile, can be introduced in extreme value index estimators. \cite{teng2021thesis} defined the empirical descending CVaR order statistics $Y_{1}^{(n)} \geqslant Y_{2}^{(n)} \geqslant \cdots \geqslant Y_{n}^{(n)}$ by
\begin{equation} \label{eqn: cvar order statistics}
    Y_{k}^{(n)}:= \frac{1}{k} \sum_{i = 1}^{k} X_{i}^{(n)}, \quad k = 1,\ldots,n,
\end{equation}
where $X_{i}^{(n)}$ are the descending order statistics. The $\{Y_{k}^{(n)} \}_{k=1}^n$
are named CVaR order statistics, and its equivalence to the empirical CVaR, i.e.,
\[
    Y_{k}^{(n)} = \widehat{\mathrm{CVaR}}_{1-\frac{k}{n}} \quad k = 1,\ldots,n.
\]
was also explained. For further relation to the CVaR and CVaR estimation we refer to \cite{rockafellar2002conditional}, \cite{acerbi2002coherence}, \cite{john2005variance}, \cite{brazauskas2008estimating}, \cite{gao2011asymptotic}, etc. Compared to VaR, CVaR captures more information about a distribution’s tail. In the field of extreme value theory, CVaR is often referred to as `conditional tail expectation' (CTE). When $X$ has a continuous distribution, CTE is a coherent risk measure and is equivalent to CVaR. Asymptotic studies of CTE in univariate or multivariate cases under related extreme-value-type conditions can be found in \cite{hua2011second}, \cite{asimit2011asymptotics}, \cite{joe2011tail}, and \cite{zhu2012asymptotic}. Instead of estimating CTE or CVaR through sample averaging, alternative estimators based on extreme value theory have been proposed for both univariate and multivariate cases. These estimators involve an intermediate step to estimate the extreme value index using methods such as the Hill estimator or the maximum likelihood approach to approximate the upper right tail of the distribution. \cite{troop2022best} introduced an estimator of CVaR in the univariate case. \cite{cai2015estimation} established an estimator of marginal expected shortfall (MES) of a random variable $X$ in a multivariate setup without imposing any parametric structure on the variables $(X,Y)$. The estimator by \cite{martins2018nonparametric} admits a nonparametric location-scale representation between two random variables. One can take advantage of extreme value theory in CVaR estimation. Conversely, one can also benefit from the CVaR quantity in extreme value index estimation. Therefore, \cite{teng2021thesis} incorporated CVaR order statistics into the framework by \cite{yun2002generalized} to introduce the CVaR-based Pickands estimator by
\begin{equation} \label{eqn: cvar estimator}
\widehat{\gamma}_{n, m}(u, v):=\frac{1}{\log v} \log \frac{Y_{m}^{(n)}-Y_{[u m]}^{(n)}}{Y_{[v m]}^{(n)}-Y_{[u v m]}^{(n)}}, \quad u, v>0, u, v \neq 1,
\end{equation}
where $1 \leqslant m,[u m],[v m],[u v m] \leqslant n$ and $[x]$ denotes the integer part of $x \in \mathbb{R}$. He showed the consistency of the estimator within $\gamma<1$. The CVaR can be viewed as an equal-weighted average of Value-at-Risk (VaR) \citep{rockafellar2002conditional}. Therefore, with $\gamma$ not too large, i.e., the distribution tail is not too heavy, the idea of equal-weighted averaging can reduce the variance of the estimator for the extreme value index (see Sections \ref{sec: large sample} and \ref{sec: simulations}).  

The rest of the paper is organized as follows. In Section \ref{sec: regular variation}, we introduce the CVaR-based regular variation condition in extreme value theory.  Section \ref{sec: cvar process} focuses on the weak convergence of empirical CVaR in functional space. Section \ref{sec: large sample} exhibits the asymptotic normality result of the CVaR-based Pickands estimator. Section \ref{sec: simulations} provides simulation studies supporting the asymptotic normality result, while Section \ref{sec: conclusion} presents the conclusions.

\section{The CVaR-based Regular Variation Conditions}
\label{sec: regular variation}
\subsection{First-order Conditions}
A sequence of positive integers $m=m(n)$ is called an intermediate sequence if $m \to \infty$ and $m / n \to 0$ as $n \to \infty$.
Define $U(x):=F^{-1}(1-1/x)$, where $F^{-1}$ denotes the quantile function of $F$. Then the condition $F\in\mathscr{D}(G_\gamma)$ is equivalent to the existence of a positive, measurable function $a(t)$ defined on a neighborhood of infinity such that
\begin{equation} \label{eqn: 1st_order_ERV_U}
    \lim_{t\to\infty} \frac{U(\frac{t}{y})-U(t)}{a(t)}=\frac{y^{-\gamma}-1}{\gamma} := h_{\gamma}(y), \quad y>0
\end{equation}
where $h_{\gamma}(y)$ has to be read as $\log{\frac{1}{y}}$ in case $\gamma=0$. Moreover, (\ref{eqn: DoA}) holds with $b_n:=U(n)$ and $a_n:=a(n)$. See Lemma 1, \cite{de1984slow} and Theorem 1.1.6, \cite{de2006extreme}. The auxiliary function $a$ in (\ref{eqn: 1st_order_ERV_U}) is regularly varying with index $\gamma$, denoted as $a\in RV_\gamma$, in other words, 
\begin{equation}\label{eqn: 1st_order_RV}
    \lim_{t\to\infty} \frac{a(\frac{t}{y})}{a(t)}=y^{-\gamma},\quad y>0.
\end{equation}
If $a\in RV_0$, $a$ is called a \textit{slowly varying} function.

Suppose $F$ has a finite mean and define 
\begin{equation} \label{eqn:V}
    V(x) := \frac{1}{x} \int_{0}^{x} F^{-1}( 1- s) ds = \frac{1}{x} \int_{0}^{x} U(\frac{1}{s}) ds.
\end{equation}
\cite{teng2021thesis} showed that if $F \in \mathscr{D}(G_\gamma)$, there exists a positive, measurable function $a(t)$ defined on a neighborhood of infinity such that
\begin{equation} \label{eqn:1st_order_ERV'_V}
\lim_{t\to\infty} \frac{V(\frac{y}{t})-V(\frac{1}{t})}{a(t)} = \frac{1}{y} \int_0^y h_\gamma(w) dw - \int_0^1 h_\gamma(w) dw := \tilde{h}_\gamma(y) \quad \text { for } y>0,
\end{equation}
alternatively, for $y>0$,
\begin{equation} \label{eqn:1st_order_ERV_V}
    \frac{V(\frac{y}{t})-V(\frac{1}{t})}{a(t)}=\tilde{h}_\gamma(y)+R(t,y),\quad R(t,y)=o(1)\:\text{as } t\to\infty,
    \tag{\ref{eqn:1st_order_ERV'_V}'}
\end{equation}
where $\tilde{h}_\gamma(y) = \frac{y^{-\gamma}-1}{\gamma(1-\gamma)}$ for $\gamma<1$ and $\tilde{h}_\gamma(y) = -log(y)$ in case $\gamma=0$.

Standing on \eqref{eqn:1st_order_ERV'_V}, \cite{teng2021thesis} proved the weak and strong consistency of the CVaR-based Pickands estimator, $\widehat{\gamma}_{n, m}(u, v)$, under the condition $F \in \mathscr{D}\left(G_{\gamma}\right)$ of \eqref{eqn: 1st_order_ERV_U} for the intermediate sequences $m=m(n)$ as $n\to\infty$ and even $m(n)/ \log (\log n) \rightarrow \infty$.

\subsection{Second-order Conditions}
To derive the asymptotic bias of the CVaR-based estimator, one needs to consider the second-order behavior of $V$ defined in \eqref{eqn:V}. The following condition first presents the second-order requirement for $U$.
\begin{condition} \label{cond1}
There are $\gamma\in\mathbb{R}, \ \rho\leq0$, and $A\in RV_{\rho}$ with $\lim_{t\to\infty}A(t)=0$ such that
\begin{equation} \label{eqn:2nd_order_ERV}
    \lim_{t\to\infty} \frac{\frac{U(\frac{t}{y}) - U(t)}{a(t)}-h_\gamma(y)}{A(t)} = H_{\gamma,\rho}(y),
\end{equation}
where the function $a$ appears in \eqref{eqn: 1st_order_ERV_U} and there exist constants $c_1$, $c_2\in\mathbb{R}$ such that
\begin{equation} \label{eqn:H}
    H_{\gamma,\rho}(y) = c_1\int_y^1 s^{-\gamma-1}\int_s^1 w^{-\rho-1} dwds + c_2 \int_y^1 s^{-(\rho+\gamma)-1} ds, \quad y>0.
\end{equation}
\end{condition}
The function $A$ appearing in \eqref{eqn:2nd_order_ERV} can be used to further restrict the sequence $m=m(n)$ and it leads to Condition 2.

\begin{condition}
\label{cond2}
Let $A$ be as in condition \ref{cond1}. The intermediate sequence $m=m(n)$ satisfies 
\[
    \sqrt{m} A(\frac{n}{m}) \to 0.
\]
\end{condition}

The second-order regular variation limit function for $V$ is similar and demonstrated in the following corollary.
\begin{corollary} \label{thm: 2nd_order_ERV_V}
If Condition \ref{cond1} holds with $\gamma<1$,
\begin{equation}
\label{eqn:2nd_order_ERV_V}  
    \lim_{t\to\infty} \frac{\frac{V(\frac{y}{t}) - V(\frac{1}{t})}{a(t)} - \tilde{h}_\gamma(y)}{A(t)}=\tilde{H}_{\gamma,\rho}(y), \quad y>0,
\end{equation}
where
\begin{equation} \label{eqn:H_tilde}
    \tilde{H}_{\gamma,\rho}(y):=\frac{1}{y}\int_{0}^{y}H_{\gamma,\rho}(x)dx - \int_{0}^{1}H_{\gamma,\rho}(x)dx.
\end{equation}
\end{corollary}

\section{The Empirical CVaR Process}
\label{sec: cvar process}

\begin{remark} \label{remark: moment}
Let $X$ be a random variable with d.f. $F$ where $F\in\mathscr{D}(G_\gamma)$, and let $k$ be an integer with $0<k<1/\gamma_+$, where $\gamma_+:=\max(0,\gamma)$. Then, $E(|X|^k)$ is finite. In particular, if $\gamma<1$, $E|X|$ is finite; if $\gamma<\frac{1}{2}$, $E(X^2)$ is finite and thus, the variance of random variable $X$ exists \citep[Theorem 5.3.1]{de2006extreme}.
\end{remark}
This implies that the CVaR-based estimators are valid for $\gamma < \frac{1}{2}$, contingent upon the finite variance of the random variable with distribution function $F$, where $F \in \mathscr{D}(G_\gamma)$, and, of course, the existence of CVaR ($\gamma < 1$).

The stochastic process with the empirical VaR process is defined by
\[
    B_n(t):=\sqrt{m}\frac{X_{[mt]}^{(n)}-U(\frac{n}{mt})}{a(\frac{n}{m})}, \quad 0<t\leq\frac{n}{m}.
\]
Assume that Condition \ref{cond1} and \ref{cond2} hold, from Proposition 1 of \cite{resnick1999smoothing}, the stochastic process
\begin{equation}
\label{eqn:process_converge}
    B_n(t):=\sqrt{m}\frac{X_{[mt]}^{(n)}-U(\frac{n}{mt})}{a(\frac{n}{m})} \Rightarrow B(t)
\end{equation}
in $D(0,\infty)$ as $n\to\infty, m\to\infty$ and $m/n\to0$, where $B(t)=t^{-\gamma-1}W(t)$ and $W$ is a standard Brownian motion.

\begin{theorem} \label{thm: CVaR process convergence}
Suppose $F\in\mathscr{D}(G_\gamma)$ with $\gamma<\frac{1}{2}$, and Conditions \ref{cond1} and \ref{cond2} hold, and define the stochastic process with empirical CVaR process by
\[
    \tilde{B}_n(t):=\sqrt{m}\dfrac{Y_{[mt]}^{(n)}-V(\dfrac{mt}{n})}{ a(\dfrac{n}{m})},\quad 0<t\leq\frac{n}{m}.
\]
Then,
\begin{equation}
\label{eqn:weak_process}
    \tilde{B}_n(t) \Rightarrow \tilde{B}(t) 
\end{equation}
in $D(0,\infty)$, where 
\begin{align*}
    \tilde{B}(t) &= \frac{1}{t} \int_0^t s^{-\gamma-1}W(s)ds, \quad t>0.
\end{align*}
\end{theorem}

Observe that when $\gamma<\frac{1}{2}$, the stochastic process $\tilde{B}(t), \ t>0$ is mean-zero normally distributed with the finite variance
\begin{equation}
\label{eqn:process_var}
    Var(\tilde{B}(t)) = \begin{cases}
        \frac{2}{(1-\gamma)(1-2\gamma)}t^{-2\gamma-1},\quad \gamma\neq0 \ \& \ \gamma<\frac{1}{2}, \\
        \frac{2}{t},\quad \gamma=0.
    \end{cases}
\end{equation}
Note that the weak convergence of the empirical CVaR holds in \eqref{eqn:weak_process} requires a finite variance of $\tilde{B}(t)$, implying $\gamma<\frac{1}{2}$. This aligns with the consistency that the CVaR-based estimators are valid when $\gamma < \frac{1}{2}$ as indicated in Remark \ref{remark: moment}. Besides, with $t_1,t_2>0$,
\begin{align*}
    Cov(\tilde{B}(t_1), \tilde{B}(t_2)) &:= \sigma(t_1,t_2) \\
    &= \begin{cases}
        (t_1 t_2)^{-1}\left[h_1(\gamma)(t_1\land t_2)^{1-2\gamma} - h_2(\gamma)(t_1\land t_2)^{1-\gamma}(t_1\lor t_2)^{-\gamma}\right], \quad \gamma\neq0 \;\&\; \gamma<\frac{1}{2} \\ \\
        (t_1\lor t_2)^{-1}\left[2 - \log{(t_1\land t_2)} + \log{(t_1\lor t_2)}\right], \quad \gamma=0,
    \end{cases}
\end{align*}
 where $h_1(\gamma)=\frac{1}{\gamma(1-\gamma)(1-2\gamma)}$ and $h_2(\gamma)=\frac{1}{\gamma(1-\gamma)}$.

\section{Asymptotic Normality of the CVaR-based Pickands Estimator}
\label{sec: large sample}
Standing on Theorem \ref{thm: CVaR process convergence} in Section \ref{sec: cvar process}, we can derive the following theorem on the asymptotic normality of the CVaR-based Pickands estimator.
\begin{theorem}[Asymptotic Normality] \label{thm: weak convergence 1}
Suppose $F\in\mathscr{D}(G_\gamma)$ with $\gamma<\frac{1}{2}$, and Conditions \ref{cond1} and \ref{cond2} hold. Let $u,v>0$ with $u,v\neq 1$ and define the function
\[
    g(v,uv) = \begin{cases}
        \frac{v^{-\gamma}-(uv)^{-\gamma}}{\gamma(1-\gamma)},\quad \gamma\neq0 \\
        \log{u}, \quad \gamma=0.
    \end{cases}
\]
Then,
\[
    \sqrt{m}(\widehat{\gamma}_{n,m}(u,v)-\gamma)\xrightarrow{d}\mathcal{N}\left(0,\frac{\sigma_\gamma^2(u,v)}{v^{2\gamma}\log^2{v}}\times \frac{1}{g^2(v,uv)}\right) \quad \text{as} \ n\to\infty,
    \qquad
\]
where
\begin{align*}
\sigma_\gamma^2(u,v) &= \begin{cases}
     \frac{2}{(1-\gamma)(1-2\gamma)}(1+v^{-1})(1+u^{-2\gamma-1})+\mathscr{G}(u,v), \quad \gamma\neq0, \\
    (1+u^{-1})(1+v^{-1}) +\mathscr{G}(u,v),\quad \gamma=0,
\end{cases}
\end{align*}
and
\[
    \mathscr{G}(u,v)=2\left(-\sigma(1,u)-v^{\gamma}\sigma(1,v)+v^{\gamma}\sigma(1,uv)+v^{\gamma}\sigma(u,v)-v^{\gamma}\sigma(u,uv)+v^{2\gamma}\sigma(v,uv)\right).
\]
\end{theorem}

\section{Simulations} \label{sec: simulations}
In the simulation study, two estimators are included: the estimator $\widehat{\gamma}_{n, m}^{Y}(u, v)$ in \eqref{eqn: yun} by \cite{yun2002generalized} and the CVaR-based Pickands estimator $\widehat{\gamma}_{n, m}(u, v)$ in \eqref{eqn: cvar estimator}. The parameters $u$ and $v$ are set to $u=v=2$, rendering the first estimator equivalent to the Pickands estimator introduced by \cite{pickands1975statistical}. The distributions involved in the experiment include Generalized Extreme Value and Generalized Pareto distributions, as they are applicable for $\gamma \in \mathbb{R}$. We compare two ratios for the two estimators. The first is the ratio of asymptotic variances, and the second is the ratio of estimated variances through simulation. For details on the asymptotic variance of the estimator $\widehat{\gamma}_{n, m}^{Y}(u, v)$, please refer to \cite{yun2002generalized}. 

\begin{figure}[H]
  \begin{minipage}{0.5\textwidth}
    \centering
    \includegraphics[width=\linewidth]{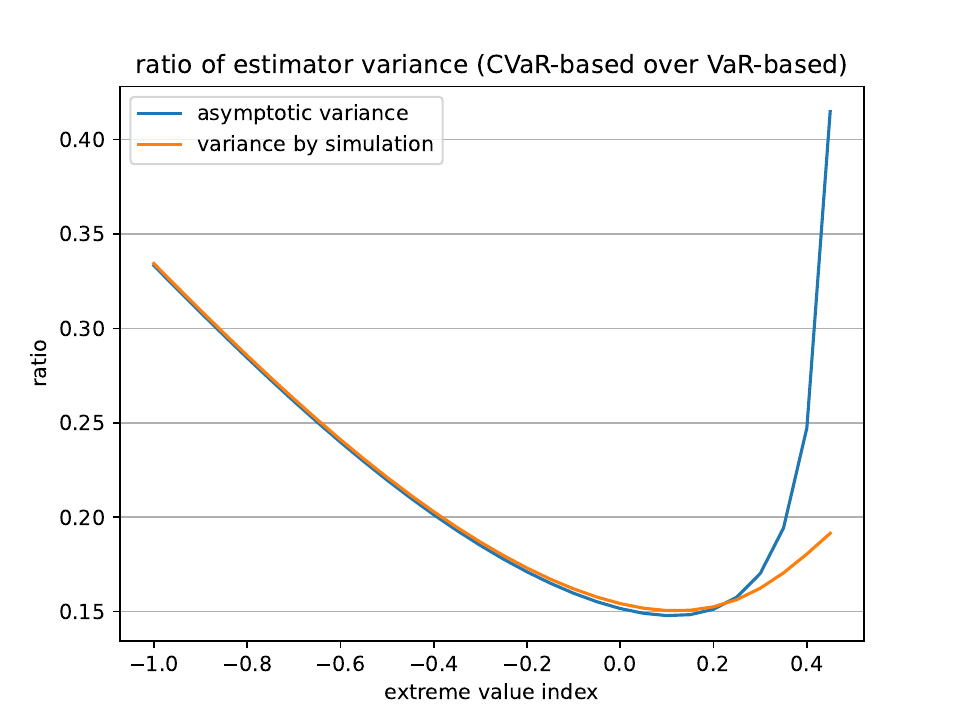}
    \caption{Generalized Extreme Value distribution.}
    \label{fig: gev ratio}
  \end{minipage}%
  \begin{minipage}{0.5\textwidth}
    \centering
    \includegraphics[width=\linewidth]{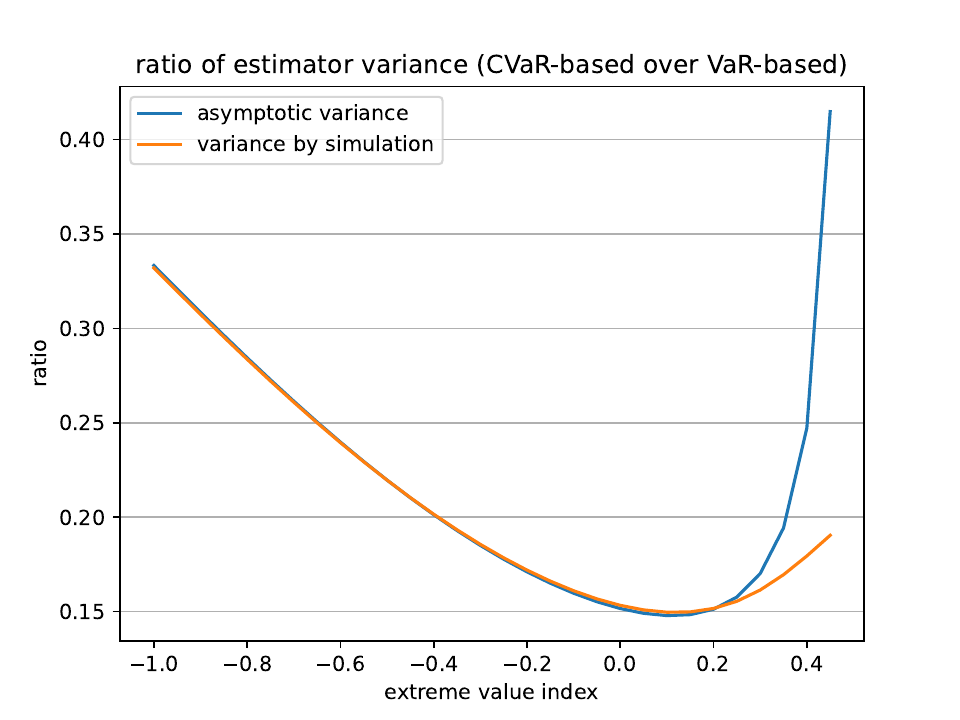}
    \caption{Generalized Pareto distribution.}
    \label{fig: gpd ratio}
  \end{minipage}
\end{figure}

The results presented in Figures \ref{fig: gev ratio} and \ref{fig: gpd ratio} are obtained from $10^4$ samples with intermediate order statistics size $m=10^2$ and sample size $n=10^4$. Each distribution is sampled from the same set of random seeds. The two figures depict the ratio of asymptotic variance and the ratio of estimated variance by simulation across the range $-1\leq\gamma<0.5$. It can be observed that the two ratios coincide within the range of $\gamma=-1$ to around $0.2$. Their discrepancies emerge and increase as $\gamma$ approaches $0.5$. This is because 
as $\gamma$ approaches $0.5$, the CVaR-based Pickands estimator eventually exhibits an infinite asymptotic variance. On the other hand, the estimated variance by simulation continues to increase but never reaches infinity. When focusing on the ratio of asymptotic variance, we can also observe that the CVaR-based Pickands estimator substantially reduces the asymptotic variance to approximately $0.15$ to $0.45$ times that of the original Pickands estimator.

\section{Conclusion} \label{sec: conclusion}
This paper studies the asymptotic behavior of the CVaR-based Pickands estimator proposed by \cite{teng2021thesis}. The CVaR-based second-order regular variation condition is derived. Then, the weak convergence of empirical CVaR in functional space is established under the derived second-order condition. Subsequently, the asymptotic normality of the estimator is presented, standing on the result of weak convergence of empirical CVaR, and is empirically validated through simulation.

\appendix \label{appendix: A}
\section{Proofs}
\noindent \textbf{Proof of Corollary \ref{thm: 2nd_order_ERV_V}}
\begin{proof}
When $\gamma<1$ and $\gamma\neq0$, since $F$ has a finite mean, both the left hand side and the right hand side in \eqref{eqn:2nd_order_ERV} are integrable functions. By Dominated Convergence Theorem, the integration is preserved by convergence. We can take integral from 0 to $y$ both left and right hand sides and take subtraction with $y=1$ to get (\ref{eqn:2nd_order_ERV_V}).
\end{proof}

\noindent \textbf{Proof of Theorem \ref{thm: CVaR process convergence}}
\begin{proof}
First note that $X_{mt}^{(n)} = F_n^{-1}(1-\frac{mt}{n})$ and $Y_{mt}^{(n)} = \frac{1}{t}\int_0^tF_n^{-1}(1-\frac{ms}{n})ds$ where $F_n^{-1}$ is the empirical quantile function for the random variable $X$. Thus, $\tilde{B}_n(t)=\frac{1}{t}\int_0^t B_n(s)ds$.

Define $\pi_t: \, D(0,\infty)\rightarrow D(0,\infty)$, where $D(0,\infty)$ is the space of cadlag functions on $(0,\infty)$, such that $\pi_t(x)=\frac{1}{t}\int_0^t x(s) ds$. The idea of the proof is to utilize the continuous mapping theorem based on the continuity of $\pi_t(x)$ on $D(0,\infty)$. On the space $D[0,k]$ with a finite and positive $k$, the Skorohod topology can be defined by the metric
\[
    d_k(x,y) = \inf_{\lambda\in\Lambda_k}{\left\{\sup_{s<t}{\left|\log{\frac{\lambda(t)-\lambda(s)}{t-s}}\right|}\lor\sup_{t}{\bigg|x(t)-y(\lambda(t))\bigg|}\right\}},
\]
where $\Lambda_k$ denotes the class of strictly increasing, continuous mappings of $[0,k]$ onto itself such that if $\lambda\in\Lambda_k$, $\lambda(0)=0$ and $\lambda(k)=k$. Meanwhile, the metric that defines the Skorohod topology on $D[0,\infty)$ is given by
\[
    d_\infty(x,y)=\sum_{k=1}^\infty 2^{-k}(1\land d_k(x^k,y^k)).
\]
From Theorem 16.1 by \cite{billingsley2013convergence}, the convergence $d_\infty(x_n,x)\to0$ in $D[0,\infty)$ if and only if there exists $\lambda_n\in\Lambda_\infty$ such that
\[
    \sup_{t<\infty}\left|\lambda_n(t)-t\right|\to0
\]
and, for each $k$,
\[
    \sup_{t\leq k}\left|(x_n \circ \lambda_n)(t)-x(t)\right| = \sup_{t\leq k}\left|x_n(\lambda_n(t))-x(t)\right|\to0,
\]
where $\Lambda_\infty$ is the set of continuous, increasing maps of $[0,\infty)$ onto itself. On $D[0,\infty)$, we define 
\[
    \pi^{'}_t(x)=\begin{cases}
        \pi_t(x), \quad t>0, \\
        0, \quad t=0.
    \end{cases}
\]
For each $k$, we have
\[
    \sup_{0<t\leq k}\left|\pi^{'}_{t}\left\{x_n \circ \lambda_n \right\}-\pi^{'}_t\left\{x\right\}\right| = \sup_{0<t\leq k}\left|\frac{1}{t}\int_0^{t}x_n(\lambda_n(s))ds - \frac{1}{t}\int_0^tx(s)ds\right|.
\]
Since any $x\in D[0,\infty)$ is bounded at any finite subinterval of $[0,\infty)$, we can find that $\sup_{0<t\leq k}\allowbreak\left|\frac{1}{t}\int_0^{t}x_n(\lambda_n(s))ds - \frac{1}{t}\int_0^tx(s)ds\right| \to 0$ follows by $\sup_{t\leq k}\left|x_n(\lambda_n(t))-x(t)\right|\to0$. When $t=0$, the convergence between $\pi^{'}_t(x_n)$ and $\pi^{'}_t(x)$ is straightforward. Thus, the convergence $d_\infty(x_n,x)\to0$ implies the convergence $d_\infty(\pi^{'}_t(x_n),\pi^{'}_t(x))\to0$ in $D[0,\infty)$. In other words, $\pi^{'}_t(x)$ is continuous on $D[0,\infty)$ and $\pi_t(x)$ is continuous on $D(0,\infty)$. Therefore, via the continuous mapping theorem (Theorem 9.7 in \citet{sen2018gentle}), we can take integral both at the LHS and RHS in (\ref{eqn:process_converge}), so we can get (\ref{eqn:weak_process}) and (\ref{eqn:process_var}). When $\gamma\neq0$, $Var(\tilde{B}(t)) = t^{-2}\int_0^t (s^{-\gamma}-t^{-\gamma})^2ds$ and the convergence of the integral requires $\gamma<\frac{1}{2}$.
\end{proof}

\noindent \textbf{Proof of Theorem \ref{thm: weak convergence 1}}
\begin{proof}
From Proposition 1, if $F\in\mathscr{D}(G_{\gamma})$ with $\gamma<\frac{1}{2}$, we have 
\begin{equation*}
\begin{aligned}
\frac{Y_{[vm]}^{(n)}-Y_{[uvm]}^{(n)}}{ a(\dfrac{n}{m})} &= \frac{V(\frac{vm}{n})-V(\frac{uvm}{n})}{ a(\dfrac{n}{m})}+\frac{1}{\sqrt{m}}(\tilde{B}_n(v)-\tilde{B}_n(uv)) \\
&= g(v,uv) + o_p(1) 
\end{aligned}
\end{equation*}
since $\tilde{B}_n(v)-\tilde{B}_n(uv)$ has asymptotically a normal distribution. Thus, as we define $A_{n,m}(u,v)=\frac{Y_{m}^{(n)}-Y_{[um]}^{(n)}}{Y_{[vm]}^{(n)}-Y_{[uvm]}^{(n)}}$,
\begin{align*}
\label{eqn:asymp}
\sqrt{m}(A_{n,m}(u,v)-v^{\gamma}) &= \sqrt{m}\frac{Y_{m}^{(n)}-Y_{[um]}^{(n)}-v^{\gamma}(Y_{[vm]}^{(n)}-Y_{[uvm]}^{(n)})}{Y_{[vm]}^{(n)}-Y_{[uvm]}^{(n)}} \\
&\distas{p} \frac{1}{g(v,uv)} \sqrt{m} \frac{Y_{m}^{(n)}-Y_{[um]}^{(n)}-v^{\gamma}(Y_{[vm]}^{(n)}-Y_{[uvm]}^{(n)})}{ a(\dfrac{n}{m})} \\
&= B_{n,m}(u,v)+C_{n,m}(u,v) \quad \text{as } n\to\infty, \numberthis
\end{align*}
where
\begin{align*}
B_{n,m}(u,v) &:= \frac{1}{g(v,uv)} \{ \tilde{B}_n(1)-\tilde{B}_n(u)-v^{\gamma}(\tilde{B}_n(v)-\tilde{B}_n(uv)) \}, \\
C_{n,m}(u,v) &:=\frac{\sqrt{m}}{g(v,uv)} \frac{V(\frac{m}{n})-V(\frac{um}{n})-v^{\gamma}(V(\frac{vm}{n})-V(\frac{uvm}{n}))}{a(\frac{n}{m})}.
\end{align*}
When $\gamma\neq0$,
\begin{align*}
    & Var\left(\tilde{B}_n(1)-\tilde{B}_n(u)-v^{\gamma}(\tilde{B}_n(v)-\tilde{B}_n(uv))\right) \\
    & \quad =\begin{pmatrix}
    1 & -1 & -v^{\gamma} & v^{\gamma}
    \end{pmatrix} \begin{pmatrix}
        h(\gamma) & \sigma(1,u) & \sigma(1,v) & \sigma(1,uv) \\
        \sigma(1,u) & u^{-2\gamma-1}h(\gamma) & \sigma(u,v) & \sigma(u,uv) \\
        \sigma(1,v) & \sigma(u,v) & v^{-2\gamma-1}h(\gamma) & \sigma(v,uv) \\
        \sigma(1,uv) & \sigma(u,uv) & \sigma(v,uv) & (uv)^{-2\gamma-1}h(\gamma)
    \end{pmatrix} \begin{pmatrix}
        1 \\ -1 \\ -v^{\gamma} \\ v^{\gamma}
    \end{pmatrix} \\
    & \quad = h(\gamma)(1+v^{-1})(1+u^{-2\gamma-1})+\mathscr{G}(u,v)
\end{align*}
and when $\gamma=0$,
\begin{align*}
    & Var\left(\tilde{B}_n(1)-\tilde{B}_n(u)-v^{\gamma}(\tilde{B}_n(v)-\tilde{B}_n(uv))\right) \\
    & \quad =\begin{pmatrix}
    1 & -1 & -v^{\gamma} & v^{\gamma}
    \end{pmatrix} \begin{pmatrix}
        2 & \sigma(1,u) & \sigma(1,v) & \sigma(1,uv) \\
        \sigma(1,u) & \frac{2}{u} & \sigma(u,v) & \sigma(u,uv) \\
        \sigma(1,v) & \sigma(u,v) & \frac{2}{v} & \sigma(v,uv) \\
        \sigma(1,uv) & \sigma(u,uv) & \sigma(v,uv) & \frac{2}{uv}
    \end{pmatrix} \begin{pmatrix}
        1 \\ -1 \\ -v^{\gamma} \\ v^{\gamma}
    \end{pmatrix} \\
    & \quad = (1+u^{-1})(1+v^{-1})+\mathscr{G}(u,v)
\end{align*}
Note here that
\begin{align*}
B_{n,m}(u,v) &\xrightarrow{d} \frac{1}{g(v,uv)} \{B(1)-B(u)-v^{\gamma}(B(v)-B(uv) \}\\
&\sim \mathcal{N}(0,\frac{\sigma_\gamma^2(u,v)}{g^2(v,uv)}) \quad \text{as } n\to\infty.
\end{align*}
Conditions (\ref{eqn:1st_order_ERV_V}) and (\ref{eqn:2nd_order_ERV_V}) imply that
\begin{equation}
\label{eqn:r_order}
C_{n,m}(u,v)=b(u, v, \gamma, \rho) \cdot \sqrt{m}A(\frac{n}{m})+\sqrt{m}o(A(\frac{n}{m})) \quad \text{as } n\to\infty,
\end{equation}
where
\[
    b(u, v, \gamma, \rho) = \begin{cases}
        \frac{\gamma(1-\gamma)}{1-u^{-\gamma}} v^{2\gamma} \left[\tilde{H}_{\gamma,\rho}(uv)-\tilde{H}_{\gamma,\rho}(v)-v^{-\gamma}\tilde{H}_{\gamma,\rho}(u) \right], \quad \gamma<1 \ \& \ \gamma\neq0, \\
        \frac{1}{\log{u}}\left[\tilde{H}_{\gamma,\rho}(uv)-\tilde{H}_{\gamma,\rho}(v)-\tilde{H}_{\gamma,\rho}(u) \right], \quad \gamma=0.
    \end{cases}
\]
From the Condition \ref{cond2}, we know that $\lim_{n\to\infty} \sqrt{m}A(\frac{n}{m})$ is equal to $0$.  Then $C_{n,m}(u,v)\to 0$ as $n\to\infty$, and so from (\ref{eqn:asymp})
\[
    \sqrt{m}(A_{n,m}(u,v)-v^{\gamma})\xrightarrow{d}\mathcal{N}(0,\frac{\sigma_\gamma^2(u,v)}{g^2(v,uv)}) \quad \text{as} \ n\to\infty. 
\] 

Define the function $h(x):=\frac{1}{\log{v}}\log{x}$, Then, applying the Taylor expansion as  $x\to v^{\gamma}$, we thus have
\begin{align*}
\sqrt{m}(\widehat{\gamma}_{n,m}(u,v)-\gamma)&=\frac{\sqrt{m}}{v^{\gamma}\log{v}}(A_{n,m}(u,v)-v^{\gamma})+\sqrt{m}o_p(A_{n,m}(u,v)-v^{\gamma}) \\
&=\frac{\sqrt{m}}{v^{\gamma}\log{v}}(A_{n,m}(u,v)-v^{\gamma})+o_p(1) \\
& \xrightarrow{d} \mathcal{N}\left(0,\frac{\sigma_\gamma^2(u,v)}{v^{2\gamma}\log^2{v}}\times \frac{1}{g^2(v,uv)}\right) \quad \text{as} \  n\to\infty.
\end{align*}
\end{proof}

\bibliography{myref}
\end{document}